\documentclass[12pt]{amsart}

\usepackage[english]{babel}
\usepackage{amsfonts}
\usepackage{amscd}
\usepackage{amsmath}

\newtheorem{rmk}{Remark}

\newtheorem{exm}{Example}

\newtheorem{prf}{Proof}

\newtheorem{prff}{}

\newtheorem{stat}{}

\newcommand{\NOT}[1]{}

\newcommand{\pa}{\par\medskip}

\newcommand{\BE}{\begin{equation}}  \newcommand{\EE}{\end{equation}}
\newcommand{\BR}{\begin{eqnarray*}}  \newcommand{\ER}{\end{eqnarray*}}
\newcommand{\BER}{$$\begin{array}}  \newcommand{\EER}{\end{array}$$}

\title{Ockham's Razor, Probability and Quantum Physics as Logic}

\author{Eliahu Levy}
\address{Department of Mathematics,
Technion -- Israel Institute of Technology,
Haifa 32000, Israel}
\email{eliahu@techunix.technion.ac.il}

\date{}

\begin{document}

\begin{abstract}
This is a philosophy-intense physics article, or, if you wish, a physics-intense philosophy article. Also, being a mathematician, I tend to view the physics, in particular the essence of quantum physics, in emphasizing the mathematical structure that serves as its language. However, I do express views on typically philosophical/epistemological matters (see, in particular, Section \ref{OckhamsRazor}). Since these points of view do not seem to me too widely expressed in the literature, while I find them quite compelling, I think this note has some interest.
\end{abstract}

\maketitle

\section{Ockham's Razor} \label{OckhamsRazor}
Let us focus on what is sometimes referred to as ``inductive logic'', i.e., the ``logic'' of knowledge, be it part of some lofty and sophisticated science, or be it an instance of the myriad manifestations of everyday knowledge. For instance, knowing that there is a table in the next room (which I do not see just now); or that our university is located in a specific town; or that it makes sense to speak of a scale of Time, described mathematically as a real line or part of it, which includes my present, past and future (totally different things from the point of view of my present consciousness); or putting other people and myself on the same footing (again, definitely not so in my ``naive'' present consciousness) and so on. \pa
It seems to me that the whole ``logic'' or ``logical process'' giving this knowledge is the principle usually referred to as Ockham's Razor. That is, a system of assertions and concepts is known, with more or less certainty, given some ``evidence'', if this system is the simplest, most harmonious, the least containing ``unnecessary complications'' of all systems that fit the evidence, the more outstanding this system is with respect to this property, the more certain it is. (Of course, that does not mean that the thus known system is, itself, simple.) \pa
The results of Ockham's Razor, thus understood, have both the flavor of ``new knowledge'', as in the examples in the first paragraph, but also the flavor of ``a description chosen to be the most economical''. Indeed, without it we would ``drown'' in possible (needlessly) complicated ``pictures'' which all fit the evidence. For example, without such ``Ockham's Razor'' we ``do not know'' that, say, all distances in the world do not really change in time in some wild manner, which we do not perceive because the measuring yardsticks change proportionally along with the measured distances. But do not mistake, we do not thus gain just ``convenient bookkeeping'' –- we gain, as in the examples in the first paragraph, what (''inductive'') knowledge is. There is no other ``reality''. (All reality being known with, possibly very high, but not absolute, certainty.)\pa
By the way, this ``inductive reasoning'' -– Ockham's Razor -- works not only with respect to the ``physical world''. One can and does apply it to abstract ``universes of discourse'' such as numbers. For example, engineers would be totally confident in using a numerical method to solve partial differential equations that came with a good rationale and had been tested on an adequate set of examples, even when no proved mathematical theorem is at hand. It is just, outstandingly, ``the simplest, most harmonious, least containing unnecessary complications, assertion'' that says that the method works. Of course, as always, that is known with some high, but not absolute, certainty. \pa
Note that the word ``harmonious'' is better not omitted here. Thus, asserting Newton's law for the force of gravity with the exponent $2$ for the distance in the denominator has a different status than with some exponent very near to $2$. If the evidence would not force us not to do so, Ockham's Razor will make us assert the exponent $2$ rather than some exponent very near to it -- the gain in doing this may be said to lie no less in harmony than in simplicity. \pa
And I stress that this same ``process'' of knowledge ``works'' both in everyday knowledge and in ``professional'' science (and, in a surprisingly similar way to many episodes of scientific discovery, in the work of real and imaginary ``Sherlock Holmses''). Science has no special method peculiar to it. Any gaining of knowledge done with supreme care (in the application of the above ``logic'', i.e.\ Ockham's Razor) deserves and is usually given the name of science. \pa

\section{Probability}
What should we make of probability (here not referring to the mathematical discipline with its characteristic language, but to its role in the sciences, in knowledge)? \pa
Some say that probability is a system of assertions of facts about the world -- say, about properties of frequencies of the different results in repeated ``experiments'' -- a branch of science on a par with, say, electrostatics. \pa
I claim that the calculus of probability does not sustain that view, as the following example suggests: \pa
Example: Toss a balanced coin independently $100$ times. The calculus of probability suggests that we will not observe a hundred times ``heads''. Rather, the frequency of ``heads'' will not deviate from $1/2$ more than a few standard deviations. Indeed, the probability for a hundred times ``heads'' is about ${10}^{-30}$. \pa
But repeat this tossing independently ${10}^{40}$ times (which is not entirely impossible physically, say in computer simulations). Then the calculus of probability predicts that in about ${10}^{10}$ of the tosses we will have a hundred ``heads''. \pa
So if probability is just some assertion of facts, it would assert simultaneously that a hundred ``heads'' must not occur and that it must occur! \pa
For me, probability is an actor on the stage of gaining knowledge; in fact, it is an instance of Ockham's Razor. \pa
Let us take an example: suppose I look from my balcony at some person passing in the street. Should I expect/assert that (s)he had been to the South Pole? If I assume that of the many identities of persons that may pass in the street below my balcony only a very small fraction had been to the South Pole, and that these identities can be considered {\em all equivalent (i.e.\ ``of the same weight'') for the consideration at hand,} then it would be an ``unnatural'', ``unharmonious'', ``involving an unnecessary complication'' assumption to say that the person I saw has ``spitefully'' fallen into the tiny fraction that had been to the South Pole! \pa
\begin{itemize}
\item  So, what probability tells us about the ``real'' world is that if we endow the space of possibilities (the sample space) with a probability measure that gives {\em equal probability to possibilities equivalent (by our considerations)} then we can assert/expect/know that an event (in the sense of probability) with very small probability will not occur. (As always, know with very high, but not absolute, certainty -- as ``inductive'' knowledge can be had at all).\pa
\item  Here, we are not just counting possibilities (putting this differently, we cannot take any probability measure for the sample space of possibilities), as a famous example shows: it is not the case that one should not expect a random student in a library to use an English book because English is just one of $6000$ existing languages!\pa
\item Similarly, the above conclusion that I should expect that the person I saw from my balcony had not been to the South Pole would not hold if I knew that there is now a parade of South Pole explorers!\pa
\item Thus, the choice of the probability measure (the decision which possibilities are ``equivalent'') is governed by our ``ambient'' knowledge/science (where we again use Ockham's Razor).\pa
\item  The particular result of the $100$ independent coin-tosses in the example above also had a probability of about ${10}^{-30}$. But it occurred.  Still, for a tossing in the future, or otherwise where we do not know the result, we should expect that this particular result shall not occur (i.e.\ that one of the other ${10}^{30}$ results will occur). Here is an instance of changing our scientific theories when new evidence makes us do so.\pa
A fact (say, the result of the above tossing) should be called {\em accidental} when we have no means to predict it in advance. (And a remark {\em in passim}: laws of nature are just known facts with some high symmetry (in space, time, etc.) or generality. Except for this symmetry/generality, they have the same status as any other knowledge. Their symmetry/generality, in itself, does not make them more ``necessary'', meaning less ``accidental'', than more particular facts.)\pa
\item To quote a famous example: suppose that a, supposedly random, picking of letters resulted in the word ``CONSTANTINOPLE''. We should weigh all our evidence to decide whether the probability measure we took was correct, and an unexpected (''accidental'') result has occurred, or we were wrong in our assessment about which possibilities were equivalent, and indeed the result should have been expected.\pa
\item  Also, what probability says about the relative frequency of occurrence of results in repeated (independent) experiments, is that we should {\em expect/assert/know} that frequency not to deviate from the probability of the result in a single experiment more than, say, a few standard deviations, since the event of it doing so has a tiny probability. (Again, as always, that knowledge is with high, not absolute, certainty.)\footnote{Defining probability as the limit of the frequency in an infinite sequence of repeated experiments (that deemed a legitimate mathematical abstraction), as sometimes proposed, does not seem to do the trick. First, that will hold only if we restrict the set of admissible infinite sequences (forbidding some set with probability zero) -- without that any limit (or divergence) may occur -- thus one here bases the {\em definition} of probability on the notion of sets of probability zero. Second, this, anyway, will not change the situation noted above, that the frequency in a finite sequence of experiments is only {\em expected} to have some relation (not to deviate more than a few standard deviations) with the probability in a single experiment.}\pa
\item  We can reason as follows about (independently) repeated experiments, say the gender of a born baby (male or female): assume we conceive a finite set $B$ of all past and future$=$all possible births, and let $100p$ be the percentage of male among them. Suppose also that the observed births are viewed as a subset of $n$ elements taken from B, and that {\em all subsets of $n$ elements of $B$ may be considered scientifically equivalent.} The set $\mathcal{E}$ consisting of all subsets of $n$ elements with percentage of male births more than a few standard deviations different from $100p$ is a tiny fraction of the set of all subsets of $n$ elements of $B$, therefore we may expect/assert/know that the percentage in our observed set does not differ too much from $100p$. Hence, we may take the observed percentage to be a suitable approximation to $100p$.\pa
\item Suppose, as with the classical case of (scientific) induction, that all observed cases have the designated property (here: being male). Then we may deduce that $p$ is very close to $1$. To put it differently: it would be an ``unnatural'' ``unharmonious'' ``involving an unnecessary complication'' assertion to say that there is a considerable possibility not to have the designated property, yet ``spitefully'' all our observed cases avoided this possibility. In fact, if there are no reasons not to say so, we will take as simplest explanation that that property is a universal law. So we will expect, with high degree of certainty, to find it also with the untested, say future, cases.\pa
    Note, however, that here we are working under the assumption that the tested and untested cases can be, to begin with, considered as ``of the same kind''/equivalent. The induction assesses this assumption along with ``inferring from the tested to the untested''. Indeed, if we begin with dissimilarity (as with Nelson Goodman's famous example of ``bleen'' defined to be blue until today and green from tomorrow onward) the inductive inference evaporates.\pa
\item One may thus answer the accusation that in applying probability to the ``real world'' we are drawing conclusions about the ``real world'' from analytic/purely mathematical statements. The non-analytical ingredient is the requirement that the probability measure in the sample space shall give equal probability to scientifically equivalent possibilities.\pa
\item In the treatment of the notion of probabilistic independence, the keen student might observe a missing ingredient. Independence is usually presented as a definition, at first sight somewhat ``arbitrary'', while on the other hand assuming/``postulating'' it, or something like it, is crucial for the presence of any ``real world'' probability-theoretical conclusions in almost all cases.\pa
The missing ingredient is that scientific/physics independence does imply probabilistic independence. If two partitions of the sample space of possibilities are scientifically independent, say pertain to two physical systems with negligible interactions, and if the elements of each partition are ``scientifically equivalent'', then so will be the elements of the partition consisting of the intersections of elements of one partition with the elements of the other. Thus in this case a product measure (in other words, assuming independence) is appropriate.\pa
\item So we have found that only events with probability very close to $0$ or very close to $1$ have any direct connection with the non-probabilistic ``world'': the former are known, with high certainty (i.e.\ as knowledge can be had), not to occur while the latter are known with high certainty to occur. We could, for that purpose, replace the probability measure we took in the sample space with, say, a measure obtained by integrating a function bounded between $1/2$ and $2$. However, to begin with we had to take a measure that gives equal weight to ``scientifically equivalent'' events, say, a measure invariant under the symmetries of the physics. In many cases this invariance fixes the measure once we know which events have small probability (recall the ergodic theorems). Also, invariant measures support the connection between the frequency of a result in repeated experiments and the probability of that result in a single experiment.\pa
But saying that an event has probability $0.7$, say, has no direct connection with the non-probabilistic world. In the actual world, the event either happens or not -- only in the case where the probability is very close to $1$ or $0$ we know, as knowledge can be had, that it will happen/will not happen. Still, saying that the probability is $0.7$ has a scientific meaning, derived from the ``scientifically'' correct probability measure on the sample space of possibilities.\pa
\end{itemize}

So, in conclusion, our vista shows an instance of Ockham's Razor -- the assertion that events with tiny probability will not occur -- obtaining a quantitative/mathematical flavor. Also, to have it we must choose a probability measure in the sample space which gives equal weight to scientifically equivalent events, say is invariant under the symmetries of the physics.

\section{Systems of Possibilities (in Classical Physics)}
It is one of the cornerstones of the way one speaks of physics to view the actual world as one ``state'' or ``possible world'' in a system of possibilities. Thus there are the laws of planetary motion under gravity, which allow many possible scenarios, one of which is the actual case. In such a ``system of possibilities'', we cannot say that logical statements are ``true'' or ``false''. Rather, they are true in some possible worlds and false in others. The logical statements (''events'' in the parlance of probability theory) form a Boolean algebra, in general bigger than just the pair $\mathbf{2}=\{\text{``true''},\text{``false''}\}$. Two statements/events correspond to the same member of the Boolean algebra if they are equivalent in the system of possibilities, i.e.\ hold or do not hold simultaneously for any state/possible world.\pa
Mathematically, we can start with the Boolean algebra, and define the states/possible worlds as Boolean homomorphisms from the Boolean algebra of events to the two-element Boolean algebra $\mathbf{2}=\{\text{``true''},\text{``false''}\}$ -- a possible world is characterized by which of the statements in the Boolean algebra hold in it. A third way will be to start with the complex ${}^*$-algebra of all bounded complex numerical magnitudes (having, in general, different numerical values in different possible worlds), with the operations of addition, multiplication, multiplication by a complex number and (complex) conjugation. Then the ``logical statements''/events will emerge as projections, i.e.\ magnitudes $p$ with the properties $p^2=p$ and $p^*=p$. The states/possible worlds will be the ${}^*$-homomorphisms (in other words, self-adjoint multiplicative linear functionals) from the algebra of bounded magnitudes to the complex numbers (mapping each magnitude to its value at that state).\pa
For Time and Space we adopt a kind of (classical!) ``Heisenberg picture''. The ``logical statements''/events and magnitudes may refer, in some manner, to a time-point, or to several time-points, and with ``naive'', i.e.\ non-general-relativistic Time, there will be the group of time-shifts acting on the logical statements, on the magnitudes, on the states/possible worlds etc., which transform each state $a$, i.e.\ each
''possible world'', to a state ``with the same relations to tomorrow as $a$ is to today''. Similarly for magnitudes etc. (And analogously there will be space-shifts by space vectors; thus the Space-Time picture of Special Relativity is readily implemented.)\pa
We may have determinism for states, i.e.\ in our system of possibilities the ``state of today determines the state of tomorrow and of all times'' (That is a ``Schr\"odinger picture''. In our ``Heisenberg picture'' we should just say that the time-shifts are uniquely defined). But the occurance of an {\em event} (in the terminology of probability), i.e.\ a set of states, equivalently a ``logical statement'', will, in general, not determine the events in other times (say the future). But if we have a probability measure (from the point of view of the algebra of bounded magnitudes -- a positive linear functional giving to $\mathbf{1}$ the value $1$) we can make statements about the probability of such future events given the present event.

\section{A Digression: Events of Probability $0$ in the Infinitary Mathematical Treatment of Probability.}
It is often essential, in science, to use what may be called ``infinitary'' mathematics even if all we ``see'' is a ``mundane'' finite system. For probability this is done by having a sample space which is a measure space of total measure $1$ with the events the $\sigma$-algebra of measurable sets. In the case of an uncountable such ``probability space'', such as the interval $[0,1]$ with Lebesgue measure, nonempty sets $E$ of measure $0$ are inevitable. Inasmuch as the points in the space are indeed considered as possible outcomes, such events $E$ are ``possible'', but in fact one does not hesitate to throw away or add such a $0$-probability event to the sample space if needed.\pa
Thus in the famous Norbert Wiener's sample space for Brownian Motion, one knows that for almost all sample points (i.e.\ for all except a set of measure $0$) the path of the motion is continuous, so one willingly takes as sample space only continuous paths.\pa
An even more striking example: in quantum mechanics, the wave-function of a particle as function of position determines the probabilities of the particle being in subsets of position space, while its Fourier transform -- the wave function of momentum -- does so for subsets of momentum space. But the Fourier transform for $L^2$ functions, such as are general wave-functions, is determined only up to a change in a set of measure $0$. Thus from the start one does not care about what happens in sets of measure $0$.\pa
In this sense, a better description would be the measure algebra or the (commutative) von~Neumann algebra $L^\infty$ of the sample space. Also there there is no restriction of countability for unions etc. -- the measure Boolean algebra is complete. The countability restriction in the usual uncountable sample space thus appears as an artifact of that way of representing the measure algebra/von~Neumann algebra, and the very status of the points of the sample space as true sample points is cast in doubt -- genuine sample points should be self-adjoint normal (i.e.\ preserving infinite monotone limits) multiplicative functionals -- states -- on the von~Neumann algebra, which are lacking in the non-discrete case.\pa
From this discussion one may conclude that is it fundamentally impossible to conceive physically an infinite set of ``occurrences'' – ``experiments'', because then any particular possible sequence of results, including the sequence occurring in the actual world, has, in general, probability $0$ and thus can be excluded from the sample space!

\section{Quantum Physics as Logic}
The advent of quantum physics came when experiments forced us to replace the commutative ${}^*$-algebra of the bounded magnitudes, as a way to define the system of possibilities, with a non-commutative ${}^*$-algebra (the algebra of the (complexified) bounded observables). This acts, above all, as a ``strange'' logic (for the system of possibilities). This logic cannot be defined, as in the commutative/classical case, by operations among the events/logical statements themselves (these being again the projections, i.e.\ members $p$ of the algebra satisfying $p^2=p$ and $p^*=p$). Indeed, Boolean operations among events=projections are defined only if they are compatible, i.e. if they commute.\footnote{One may be tempted to extend the definition of union and intersection to non-compatible events as the sum and intersection of the relevant subspaces of the Hilbert space. Note, however, that these operations depend highly non-continuousely on the subspaces, which seems to positively disqualify them.}
 The logic is rather given by the non-commutative ${}^*$-algebra itself, whose operations are always defined. (Thus there is no harm if for defining the logic we take the algebra of {\em bounded} observables, while observables in general need not be bounded.) To decipher what this logic says we must proceed, as much as possible, by analogy with ``usual'' (i.e.\ commutative algebra) logic.\pa
Note, that this is just the logic that science teaches us to put in the system of possibilities. It is not a rival to the usual logic of, say, mathematics.  Indeed, in treating that system mathematically we are with ordinary mathematical logic.\pa
In the commutative case probability measures were positive linear functionals\footnote{In fact, to be mathematically correct, in the infinite-dimensional case, only positive linear functionals belonging to the {\em predual} (of the von Neumann algebra) and mapping $\mathbf{1}$ to $1$ are to be taken as probability measures. The same in the non-commutative case.} on the algebra mapping $\mathbf{1}$ to $1$. By analogy, here also we view positive linear functionals mapping $\mathbf{1}$ to $1$ as ``probability measures''. Such a probability measure gives to each event a probability $0\le p\le1$, which for disjoint compatible events (for whom a union is defined) will be additive. Of course, when the non-commutative algebra is written as an algebra of operators in a Hilbert space, these ``probability measures'' take the form of density matrices (positive Hermitian operators of trace $1$).\pa
The non-commutativity of the algebra changes completely the role of states. Genuine states, following the commutative case, should be ${}^*$-homomorphisms from the algebra of (the complexified) bounded observables to the complex numbers (i.e.\ multiplicative self-adjoint linear functionals). There every observable will get its value and every projection=event will get the value $1$ or $0$ $=$ ``true'' or ``false''. But for a non-commutative algebra such ${}^*$-homomorphisms are rare and insufficient. The ``states'' spoken of in quantum physics, given by vectors in a Hilbert space up to multiplication by a scalar, are something else. These are events characterized by being ``atomic'' in a mathematical sense, i.e.\ they have no proper sub-events. Equivalently, there is only one ``probability measure'' supported in them. Physically, they give the maximum specificity that one can have, what in the classical/commutative case had characterized single states/possible worlds. But their unique probability measure gives to a general event a value different from $0$ or $1$. In particular, future events have just probabilities with respect to such ``mathematically atomic'' ``states''. In this these maximum specificity states are like general events -- sets of states -- in the commutative case.\pa
Note, that in the commutative case only states/possible worlds could endow statements/events of the system of possibilities with a truth-value 'true'' or ``false'', so that these events occur or not/these statements are true or false. In the system of possibilities itself these are just elements of a Boolean algebra, and saying that they occur or not is meaningless. The same holds for the non-commutative quantum system of possibilities.\pa
But, in this quantum picture, we still have to recover our ``actual world'' -- where events occur or not -- and our usual Boolean logic. It seems clear where to find them. We need a commutative sub-algebra, and something like that presents itself: the algebra of the macroscopic, quasi-classical observables which {\em almost} commute. These, and the events they define, are what we have in our old classically behaving world. Note that the evolution of this system in Time is defined by conjugation with imaginary exponents of the energy (Schr\"odinger's equation), So the energy, itself quasi-classical, cannot exactly commute with other quasi-classical observables, otherwise there would be no time evolution there. Similarly with the momentum which induces variation in space. Thus this quasi-classical sub-algebra is {\em approximately} commutative, (in fact, even its members are only approximate -- they cannot be handled in greater precision than the ``uncertainty'' that makes them commute) and our ``actual world'' is described by an approximate *-homomorphism from this algebra to the complex numbers, which will give values $0$ or $1$ to projections (events) belonging to this ``quasi-classical'' algebra, i.e.\ truth values ``true'' or ``false'' to these ``quasi-classical'' statements. Only these events occur or not in our ``actual'' world, and only with them we can use our usual logic. Their approximate nature is usually unnoticed by us, since we ourselves come from this ``approximate world'', but it limits the number of different events (statements) that we can meaningfully conjunct or disjunct, thus limits the number of things -- amount of information -- that we can speak about (to something like a ``mundane'' action measured in Planck's Constant -- something like ${10}^{34}$), limits the amount of time to the past or future that can have meaning for our ``actual world'' (because the non-commutativity with observables transformed by Hamiltonian evolution, although small and negligible for our mundane intervals of Time, becomes big for ``enormous'' intervals, hence one cannot include presumably ``quasi-classical'' quantities pertaining to enormously distant times or distances -- something like multiples of mundane times or distances by the ratio of a mundane action to Planck's constant -- in the same approximate commutative algebra), all that making our physical ``actual world'' finite to a delimited extent. The infinite space or time models used in physics serve, in this respect, a similar role to the infinite plane of coordinates in which a map includes the grounds of a city.\pa
And a great origin of ``paradoxes'' in thinking about Quantum Theory is our reluctance to obey its ``strange'' logic, where it is meaningful for events to happen or not only when speaking about the ``actual world'', that ``extra ingredient'' described by an approximate ${}^*$-homomorphism from the approximately commutative algebra of the quasi-classical observables to the complex numbers, while without that extra ingredient we have the Theory with its non- commutative logic, dealing only with the ``system of possibilities'', where it is meaningless for events to occur or not. But this Theory does speak about particles, fields, physical systems etc.\ and one is so tempted to say, in its frame, that ``the electron {\em is} here or there'' ``it {\em has} this or that property'' ``the system {\em is} (or was) in that state'' as if events there happened or not, which the ``non-commutative logic'' of the ``system of possibilities'' forbids. And then one runs straight into paradoxes.\pa
Note that it seems that we really have here an ``extra ingredient''. One might wonder whether we could not deduce everything in our ``actual world'' from ``probability close to $1$'' arguments. This seems not to be the case. It seems that in many cases quantum fluctuations have been magnified to macroscopic consequences, making many different outcomes each with small probability, of which just one is asserted in the actual world.  And moreover there are so many details in our ``actual world'' that seem entirely erratic.\pa
Of course, we can investigate non-quasi-classical systems only by making them bear on our almost-commutative quasi-classical ``world'', i.e.\ by measuring them. Moreover, our ``world'' is protected from ``stray non-commutativity'', such as carrying conclusions of former measurements to the future via the Hamiltonian (Schr\"odinger's equation) evolution, by decoherence, which will wipe out any such conclusions, preserve only what is quasi-classical and (almost) commuting and thus create the separating wall between the quasi-classical and the truly quantum worlds. Any ``Schr\"odinger cat'' (or ``Schr\"odinger physicist or mathematician'', for that matter) is either in the quasi-classical domain, hence one may assume in principle that in our ``actual world'' the question: is (s)he alive? is settled, or is in the truly quantum domain, where superpositions are routine, but can then be investigated by us only via measurements.\pa
Usually, the quasi-classical world is governed by the deterministic laws of classical physics, to be derived, in principle, from the quantum theory. But the fact that everything is approximate has consequences. Thus when the deterministic classical equations are chaotic we have true non-determinism in the quasi-classical system: between assertions about far enough time-moments one may have only probabilistic relations. Another case of non-determinism comes from measurements and measurements-like phenomena, where ``truly quantum'' elements bear on the quasi-classical ``world''.\pa
When we proceed to apply Ockham's Razor, in its particular ``probabilistic'' flavor, to such quantum systems, we have to follow the ``non-commutative'' guidelines. We will have a ``probability measure'', dictated by the science and its symmetries, in the non-commutative system of possibilities, i.e.\ a positive linear functional on the non-commutative algebra of (complexified) bounded observables giving to $\mathbf{1}$ the value $1$ ($=$ a density matrix). Its restriction to the approximately-commutative quasi-classical algebra will give an ordinary (approximate) probability measure on the quasi-classical (approximate) events, and if some (quasi-classical) event has tiny conditional probability relative to the (quasi-classical) things that we already know to be true in the ``actual world'', then we assert/know (with high certainty) that it will not occur (in our actual world). Indeed, if it had exactly zero conditional probability, then the theory would say that its negation includes = follows from what we already know. Here, since the event has tiny conditional probability, we should not assume that our actual world has ``spitefully'' fallen into it.\pa
(Note that here we work in the general framework of quantum physics, the ``acceptance'' of which, of course, is also a result of ``inductive logic'' $=$ Ockham's Razor, which ``produces'' all our (''inductive'', known with high certainty) knowledge. That we must, in principle, be prepared to ``amend'' if forced to do so by further future evidence.)\pa
Let us consider the proverbial quantum measurement. A quantum system is prepared, by making a measurement and taking only the cases with suitable outcomes (say, ``electrons that move in a certain way''), Then, maybe after some ``development'' in the quantum system, another measurement is made, then, maybe, more measurements. We, of course, ``live'' in the quasi-classical world, where in our ``actual world'' the measuring apparatuses recorded results of the measurements. Here, contrary to a quantum ``non-commutative system of possibilities'', assertions are true or false -- it is true that the measurements gave these results. Note that the total Hilbert space must have room for keeping records of all these measurement results. In particular, it must keep a record that the system had been appropriately prepared (by a measurement and selecting the desirable cases).\pa
We assume that each measurement distinguished between all vectors in a basis of the Hilbert space of the quantum system to be measured (there is no ``classical randomness'').\pa
The computations of what should be expected in this experiment would follow, of course, the usual rules of Copenhagen Quantum Theory. This is because mathematically, each subsequent measurement, by adding to the quasi-classical record, basically makes our entire Hilbert space tensor with the Hilbert space of the new (quasi-classical) record. In writing the ``probability distribution'' -- density matrix -- that we should take, always in a basis that diagonalizes the quasi-classical observables, each entry $w_{ij}$ of the density matrix will be replaced, by the above tensoring, with a sub-matrix with trace $w_{ij}$, (a positive Hermitian sub-matrix in the case of a diagonal element $w_{ii}$), and conditioning to the quasi-classical event of the specific result of the new measurement (true in the actual world) picks a particular entry in the diagonal of that sub-matrix. For the next measurement this entry itself is to be expanded into a sub-matrix, etc. Note the remarkable fact that every subsequent measurement as if ``imposes'' its pure state as the new density matrix after the conditioning -- we can ``forget'' what density matrix we used before. In other words, we must compute as if each measurement had induced a collapse of the state of the measured quantum system (which has nothing to do with the group of time-shifts -- Hamiltonian evolution -- which of course, always acts by Schr\"odinger's equation).\pa
So what does Ockham's Razor tell us here?\pa
To speak loosely, the result of every measurement has singled out a pure state $\tau$ which supports only one ``probability distribution'' -- density matrix, i.e.\ $\tau^*\otimes\tau$. Hence we know ``which probability distribution to take''. And we find that for independent repetitions of the experiments, the quasi-classical event that indicates that the correct frequency of a property $=$ measurement result (in the correct margin of error) will occur, is ``almost'' implied by $=$ contains the future of the event stating that we prepared the experiments correctly -- {\em ``almost'' -- in the sense of the probability distribution dictated by the ``pure states'' of the measurement results, which is, indeed, the probability measure that we should take}. Hence we must expect that this event will occur (in our actual world) -- our actual world will not ``spitefully'' fall into the ``negligible by weight'' negation of that event -- while we do not have such ``probability close to $0$ or $1$'' situation for the result of each particular experiment, so there we do not ``know''. That is, of course, the usual scenario, also in classical physics, with probabilistic, non deterministic cases.

\end{document}